\documentclass[12pt]{article} 
\usepackage{amsmath,amssymb}
\usepackage{rotating}
\usepackage{amssymb}
\usepackage{amsmath,amscd}
\usepackage{pst-node,pstricks,multido,pst-plot,pst-text,pst-3d}%
\usepackage{graphicx}
\usepackage{amsfonts}
\usepackage{url}

\usepackage[latin1]{inputenc}
\usepackage{subfigure}
\usepackage[round,numbers]{natbib}
\newtheorem{example}{Example}[section]

\newtheorem{theorem}[example]{Theorem}

\newtheorem{definition}[example]{Definition}
\newtheorem{proposition}[example]{Proposition}

\def\S{{\mathfrak  S}}

\def\<{\langle}
\def\>{\rangle}

\def\C{{\mathbb C}}

\def\N{{\mathbb N}}

\def\ie{{\it i.e.}, }

\def\goth{\mathfrak}

\def\shuff#1#2{\mathbin{
\hbox{\vbox{ \hbox{\vrule \hskip#2 \vrule height#1 width 0pt
}%
\hrule}%
\vbox{ \hbox{\vrule \hskip#2 \vrule height#1 width 0pt
\vrule }%
\hrule}%
}}}

\def\shuf{{\mathchoice{\shuff{7pt}{3.5pt}}%
{\shuff{6pt}{3pt}}%
{\shuff{4pt}{2pt}}%
{\shuff{3pt}{1.5pt}}}}%
\def\shuffle{\,\shuf\,}

\def\ashuff#1#2#3{
\kern 1pt \vrule height#1 \overline{\vrule height#3 width 0pt
\hskip#2} \rule{.3pt}{#1}\overline{\vrule height#3 width 0pt
\hskip#2} \rule{.3pt}{#1} \kern 1pt }

\def\Carre3#1{\left[\begin{array}{ccc}#1\end{array}\right]}

\def\S{{\goth S}}
\def\WSym{\mathrm{WSym}}
\def\nul{{\goth o}}
\def\type{\mathrm{type}}
\def\Stab{\mathrm{Stab}}
\def\orb{\mathrm{orb}}
\def\BWSym{\mathrm{BWSym}}
\def\FQSym{\mathrm{FQSym}}

\def\BSym{\mathrm{BSym}}
\def\BQSym{\mathrm{BQSym}}
\def\BWQSym{\mathrm{BWQSym}}
\def\Sym{\mathrm{Sym}}
\def\QSym{\mathrm{QSym}}
\def\WQSym{\mathrm{WQSym}}
\def\std{\mathrm{std}}

\begin{document}
\title{Word symmetric functions and the Redfield-P\'olya  theorem\thanks{This paper is partially supported by the ANR project PhysComb, ANR-08-BLAN-0243-04, the CMEP-TASSILI project 09MDU765 and the ANR project CARMA ANR-12-BS01-0017.}}
\author{\small Jean-Paul Bultel \and  \small Ali Chouria \and \small Jean-Gabriel Luque \and \small Olivier Mallet\\ \small Laboratoire LITIS - EA 4108 Universit\'e de Rouen,\\ \small Avenue de l'Universit\'e,
 76801 Saint-\'Etienne-du-Rouvray Cedex.}
\maketitle
\begin{abstract}
We give noncommutative versions of the Redfield-P\'olya theorem in $\mathrm{WSym}$,
 the algebra of word symmetric functions, and in other related combinatorial Hopf algebras.\\
 \end{abstract}

\section{Introduction}

The Redfield-P\'olya enumeration theorem (R-P theorem) is one of the most exciting  results in combinatorics of the twentieth century. It was first published by John Howard Redfield \cite{Redfield} in 1927 and independently rediscovered by George P\'olya ten years later \cite{Polya}. Their motivation was to generalize Burnside's lemma on the number of orbits of a group action on a set (see \emph{e.g} \cite{Burnside}). Note that Burnside attributed this result to Frobenius \cite{Frobenius} and  it seems that the formula was prior known to Cauchy. Although Redfield found the theorem before P\'olya, it is often attributed only to P\'olya. This is certainly due to the fact that P\'olya popularized the result by providing numerous applications to counting problems and in particular to the enumeration  of chemical compounds. 
The theorem is a result of group theory but there are important implications in many disciplines (chemistry, theoretical physics, mathematics --- in particular combinatorics and enumeration \emph{etc.}) and its extensions lead to Andr\'es Joyal's combinatorial species theory \cite{BLL}.\\
Consider two sets $X$ and $Y$ ($X$ finite) and let $G$ be a finite group acting on $X$. For a map $f:X\rightarrow Y$, define the vector $v_f=(\#\{x:f(x)=y)\})_{y\in Y}\in \N^Y$. The R-P theorem deals with the enumeration of the maps $f$ having a given $v_f=v$ ($v$ fixed) up to the action of the group $G$. The reader can refer  to \cite{Krishna} for proof, details, examples and generalizations of the R-P theorem.\\
Algebraically, the theorem can be pleasantly stated in terms of symmetric functions (see \emph{e.g.} \cite{Macdo,Lascoux}); the cycle index polynomial is defined in terms of power sum symmetric functions and, whence writing in the monomial basis, the coefficients count the number of orbits of a given type. From the seminal paper \cite{NCSFI}, many combinatorial Hopf algebras have been discovered and investigated. The goal is to mimic the combinatorics and representation theory related to symmetric functions in other contexts.
This paper asks the question of the existence of combinatorial Hopf algebras in which the R-P theorem can naturally arise. In this sense, the article is the continuation of \cite{DLPT,DLNTT} in which the authors investigated some Hopf algebras in the aim to study the enumeration of bipartite graphs up to the permutations of the vertices.\\
Each function $f:X\rightarrow Y$ can be encoded by a word of size $\#X$ on an alphabet $A_Y=\{a_y:y\in Y\}$. Hence, intuitively, we guess that the R-P theorem arises in a natural way in the algebra of words $\C\langle A\rangle$. The Hopf algebra of word symmetric functions $\WSym$ has been studied in \cite{RS,BRRZ,HNT}. In S
Section \ref{SWSym}, we recall the basic definitions and properties related to this algebra and propose a definition for the specialization of an alphabet using the concept of operad \cite{LV,Markl}. In Section \ref{BSym}, we construct and study other related combinatorial Hopf algebras. In Section \ref{Polyaetc}, we investigate the analogues of the cycle index polynomials in these algebras and give two noncommutative versions of the R-P theorem. In particular, we give a word version and a noncommutative version. 
Finally, in Subsection \ref{Harary}, we propose a way to raise Harary-Palmer type enumerations (the functions are now enumerated up to an action of $G$ on $X$ and an action of another group $H$ on $Y$) in $\WSym$. For this last equality, we need the notion of specialization defined in Section \ref{SWSym}.

\section{Word symmetric functions\label{SWSym}}

\subsection{Basic definitions and properties}


Consider the family $\Phi:=\{\Phi^\pi\}_\pi$ whose elements are indexed by set partitions of $\{1,\ldots,n\}$ (we will denote $\pi\Vdash n$).
 The algebra $\WSym$ \cite{RS} is generated by $\Phi$ for the shifted concatenation product:
$
\Phi^\pi\Phi^{\pi'} = \Phi^{\pi\pi'[n]}
$
where $\pi$ and $\pi'$ are set partitions of $\{1,\ldots,n\}$ and $\{1,\ldots,m\}$, respectively, and $\pi'[n]$ means that we add $n$ to each integer occurring in $\pi'$.
Other bases are known, for example, the word monomial functions defined by
\[
\Phi^\pi = \sum_{\pi \le \pi'} M_{\pi'}
\]
where $\pi \le \pi'$ indicates that $\pi$ is finer than $\pi'$, \ie that each block of $\pi'$ is a union of blocks of $\pi$.

$\WSym$ is a Hopf algebra when endowed with the shifted concatenation product and the following coproduct, where $\mathrm{std}(\pi)$ means that for all $i$, we replace the $i$th smallest integer in $\pi$ by $i$:
\[
\Delta M_\pi = \sum_{\substack{\pi'\cup\pi''=\pi\\\pi'\cap\pi''=\varnothing}} M_{\mathrm{std}(\pi')} \otimes M_{\mathrm{std}(\pi'')}.
\]
Note that the notion of standardization makes sense for more general objects. 
If $S$ is a total ordered set, the standardized $\std(\ell)$, 
for any list $\ell$ of $n$ elements of $S$, is classicaly the permutation $\sigma=[\sigma_1,\dots,\sigma_n]$  
verfying $\sigma[i]>\sigma[j]$ if $\ell[i]>\ell[j]$ or if $\ell[i]=\ell[j]$ and $i>j$. Now if the description of an object $o$
 contains a list $\ell$ , the standardized $\std(o)$ is obtained by replacing $\ell$ by $\std(\ell)$ in $o$.

Let $ A $ be an infinite alphabet. The algebra $\WSym$ is isomorphic to $\WSym( A )$, 
the subalgebra of $\C\langle  A \rangle$ defined by Rosas and Sagan \cite{RS} and constituted by the polynomials which are invariant by permutation of the letters of the alphabet. The explicit isomorphism sends each $\Phi^\pi$ to the polynomial
$
\Phi^\pi( A ) := \sum_w w
$
where the sum is over the words $w = w_1\cdots w_n$ ($w_1,\ldots,w_n \in A$) such that if $i$ and $j$ are in the same block of $\pi$ then $w_i$ = $w_j$.
Under this isomorphism, each $M_\pi$ is sent to
$
M_\pi ( A ) = \sum_w w
$
where the sum is over the words $w = w_1\cdots w_n$ ($w_1,\ldots,w_n \in  A $)
 such that $w_i$ = $w_j$ if and only if $i$ and $j$ are in the same block of $\pi$.
In the sequel, when there is no ambiguity, we will identify the algebras $\WSym$ and $\WSym( A )$.
With the realization explained above, 
the coproduct of $\WSym$ consists of identifying the algebra 
$\WSym\otimes\WSym$ with $\WSym( A + B )$, where 
$ A $ and $ B $ are two  noncommutative alphabets such that $ A $ commutes with $ B $,
 by setting $f( A )g( B )\sim f\otimes g$ (see. 
It is a cocommutative coproduct for which the polynomials $\Phi^{\{1,...,n\}}$ are primitive. Endowed with this coproduct, $\WSym$ has a Hopf structure which has been studied by Hivert \emph{et al.} \cite{HNT} and Bergeron \emph{et al.} \cite{BRRZ}.

\subsection{What are virtual alphabets in $\WSym$?}

We consider the set $\goth C$ of set compositions together with additional elements $\{\nul_m:m>1\}$ 
(we will also set $\nul_0=[]$) and a unity ${\mathbf 1}$. 
This set is a naturally bigraded set: 
if ${\goth C}^m_n$ denotes the set of compositions of $\{1,\dots,n\}$ into $m$ subsets, we have
$
{\goth C}=\{\mathbf 1\}\cup\bigcup_{n,m\in\N}{\goth C}^m_n\cup\{\nul_m:m>1\}.$
We will also use the notations ${\goth C}_n$ (resp. ${\goth C}^m$) to denote the set of compositions of $\{1,\dots,n\}$ (resp. the set compositions into $m$ subsets together with $\nul_m$) with the special case:
$
\mathbf 1\in{\goth C}^1.
$
The formal space $\C[{\goth C}^m]$ is naturally endowed with a structure of right $\C[\S_m]$-module; the permutations acting by permuting the blocks of each composition and letting $\nul_m$ invariant. For simplicity, we will denote also by $\goth C$ the collection ($\mathbb S$-module, see \emph{e.g.}\cite{LV})$ [\C[{\goth C}^0],\C[{\goth C}^1],\dots,\C[{\goth C}^m]]$.

For each $1\leq i\leq k$, we define partial compositions $\circ_i:{\goth C}^k\times{\goth C}^{k'}\rightarrow {\goth C}^{k+k'-1}$ by:
\begin{enumerate}
\item If $\Pi=[\pi_1,\dots,\pi_k]$ and $\Pi'=[\pi'_1,\dots,\pi'_{k'}]$ then
\[
 \Pi\circ_i \Pi'=\left\{
\begin{array}{ll}
[\pi_1,\dots,\pi_{i-1},\pi'_1[\pi_{i}],\dots,\pi'_{k'}[\pi_i],\pi_{i+1},\dots,\pi_k]&\mbox{ if }
\Pi'\in {\goth C}_{\#\pi_i}\\
\nul_{k+k'-1}&\mbox{ otherwise},
\end{array}\right.
\]
where $\pi'_{j}[\pi_i]=\{i_{j_1},\dots,i_{j_p}\}$ if $\pi'_j=\{j_1,\dots,j_p\}$ and $\pi_i=\{i_1,\dots,i_k\}$ with $i_1<\dots<i_k$;
\item $\Pi\circ_i \nul_{k'}=\nul_k\circ_i\Pi'=\nul_{k+k'-1}$ for each $\Pi\in{\goth C}^k$ and $\Pi'\in{\goth C}^{k'}$;
\item ${\mathbf 1}\circ_1 \Pi'=\Pi'$ and $\Pi\circ_i{\mathbf 1}=\Pi$ for each $\Pi\in{\goth C}^k$ and $\Pi'\in{\goth C}^{k'}$.
\end{enumerate}
\begin{proposition}
The $\mathbb S$-module $\goth C$ (\emph{i.e.} each graded component ${\goth C}_n$ is a $\S_n$-module \cite{LV}) endowed with the partial compositions $\circ_i$ is an operad in the sense of Martin Markl \cite{Markl}, which means the compositions satisfy:
\begin{enumerate}
\item (Associativity) For each $1\leq j\leq k$, $\Pi\in{\goth C}^k$, $\Pi'\in{\goth C}^{k'}$ and $\Pi''\in{\goth C}^{k''}$:
\[
(\Pi\circ_j \Pi')\circ_i\Pi''=\left\{
\begin{array}{ll}
	(\Pi\circ_i\Pi'')\circ_{j+k''-1}\Pi',&\mbox{ for }1\leq i< j,\\
\Pi\circ_J(\Pi'\circ_{i-j+1}\Pi''),&\mbox{ for }j\leq i<k'+j,\\
(\Pi\circ_{i-k'+1}\Pi'')\circ_{j}\Pi',&\mbox{ for }j+k'\leq i< k+k'-1,\\
\end{array}
\right.
\]
\item (Equivariance) For each $1\leq i\leq m$,  $\tau\in\S_m$ and $\sigma\in S_n$, let $\tau\circ_i\sigma\in \S_{m+n-1}$ be given by
inserting the permutation $\sigma$ at the $i$th place in $\tau$. 
If $\Pi\in{\goth C}^k$ and $\Pi'\in{\goth C}^{k'}$ then $(\Pi\tau)\circ_i(\Pi'\sigma)=(\Pi\circ_{\tau(i)}\Pi')(\tau\circ_i\sigma)$.
\item (Unitality) ${\mathbf 1}\circ_1 \Pi'=\Pi'$ and $\Pi\circ_i{\mathbf 1}=\Pi$ for each $\Pi\in{\goth C}^k$ and $\Pi'\in{\goth C}^{k'}$.
\end{enumerate}
\end{proposition}
Let $V=\bigoplus_n V_n$ be a graded space over $\C$. We will say that $V$ is a (symmetric) $\goth C$-module, if there is an action of $\goth C$ on $V$ which satisfies:
\begin{enumerate}
\item Each $\Pi\in\goth C^m$ acts as a linear application $V^m\rightarrow V$. 
\item (Compatibility with the graduation) 
$[\pi_1,\dots,\pi_m](V_{j_1},\dots,V_{j_m})={\bf 0}$ if $\#\pi_i\neq j_i$ 
for some $1\leq i\leq m$ 
and $\nul_m(V_{j_1},\dots,V_{j_m})={\bf 0}$. 
Otherwise, $[\pi_1,\dots,\pi_m]\in{\goth C}_n$ 
sends $V_{\#\pi_1}\times\dots\times V_{\#\pi_m}$ to $V_n$. Note also the special case $\mathbf 1(v)=v$ for each $v\in V$.
\item (Compatibility with the compositions)
$$(\Pi\circ_i\Pi')(v_1,\dots,v_{k+k'-1})=\Pi(v_1,\dots,v_{i-1},\Pi'(v_i,\dots,v_{i+k'-1}),
v_{i+k'},\dots,v_{k+k'-1}).$$
\item\label{sym} (Symmetry) $\Pi(v_1,\dots,v_k)=(\Pi\sigma)(v_{\sigma(1)},\dots,v_{\sigma(k)})$.
\end{enumerate}
If  $V$ is generated (as a $\goth C$-module) by \ \  $\{v_n:n\geq1\}$ with $v_n\in V_n$, then setting $v^{\{\pi_1,\dots,\pi_m\}}=[\pi_1,\dots,\pi_m](v_{\#\pi_1},\dots,v_{\#\pi_m})$, we have $V_n=\mathrm{span}\{v^\pi:\pi\Vdash n\}$. Note that the existence of $v^\pi$ follows from the point \ref{sym} of the definition of a $\goth C$-module.
\begin{example}\label{refshuffle}\rm
If $A$ is a noncommutative alphabet, the algebra $\C\langle A\rangle$ can be endowed with a structure of $\goth C$-module by setting
\[
[\pi_1,\dots,\pi_m](w_1,\dots,w_m)=\left\{\begin{array}{ll}
\shuffle_{[\pi_1,\dots,\pi_m]}(w_1,\dots,w_m)&\mbox{ if } w_i\in A^{\#\pi_i}\mbox{ for each }1\leq i\leq m\\
0&\mbox{ otherwise}\end{array}
\right.
\] 
where $w_1,\dots, w_m\in A^*$ and $\shuffle_{[\pi_1,\dots,\pi_m]}(w_1,\dots,w_m)=a_1\dots a_n$ is the only word of $A^{\#\pi_1+\dots+\#\pi_m}$ such that for each $1\leq i\leq m$, if $\pi_i=\{j_1,\dots,j_\ell\}$, $a_{j_1}\dots a_{j_\ell}= w_j$.\\
Note that $\C[A]$ has also a structure of $\goth C$-module defined by $[\pi_1,\dots,\pi_m](x_1,\dots,x_m)=x_1\dots x_m$ if $x_i$ is a monomial of degree $\#\pi_i$.
\end{example}
\begin{proposition}
$\WSym$ is a ${\goth C}$-module.
\end{proposition}
{\bf Proof:}\\
We define the action of ${\goth C}$ on the power sums by 
\[
[\pi_1,\dots,\pi_m](\Phi^{n_1},\dots,\Phi^{n_m})=
\left\{\begin{array}{ll}
\Phi^{\{\pi_1,\dots,\pi_m\}}&\mbox{ if } n_i=\#\pi_i\mbox{ for each }1\leq i\leq m\\
0&\mbox{ otherwise}.
\end{array}\right.
\]
and extend it linearly to the spaces $\WSym_n$. Since this action is compatible with the realization:
\[
[\pi_1,\dots,\pi_m](\Phi^{n_1},\dots,\Phi^{n_m})(A)=\shuffle_{[\pi_1,\dots,\pi_m]}(\Phi^{n_1}(A),\dots,\Phi^{n_m}(A))
\]
(the definition of $\shuffle_\pi$ is given in Example \ref{refshuffle}) and $\WSym$ is obviously stable by the action of $\goth C$, $\WSym(A)$ is a sub-$\goth C$-module of $\C\langle A\rangle$. Hence, $\WSym$ is a $\goth C$-module.
$\Box$\\ \\
A morphism of $\goth C$-module  is a linear map $\varphi$ from a $\goth C$-module $V_1$ to another $\goth C$-module $V_2$ satisfying $\varphi\Pi(v_1,\dots,v_m)=\Pi(\varphi(v_1),\dots,\varphi(v_m))$ for each $\Pi\in{\goth C}^m$.
In this context, a virtual alphabet (or a specialization) is defined by a morphism of $\goth C$-module  $\varphi$ from $\WSym$ to a $\goth C$-module $V$. The image of a word symmetric function $f$ will be denoted by $f[\varphi]$. 
The $\goth C$-module $\WSym$ is free in the following sense:
\begin{proposition}
Let $V$ be  generated by $v_1\in V_1, v_2\in V_2,\dots, v_n\in V_n, \dots$ as a $\goth C$-module (or equivalently, $V_n=\mathrm{span}\{v^\pi:\pi\Vdash n\}$). There exists a morphism of $\goth C$-module $\varphi: \WSym\rightarrow V$, which sends $\Phi^n$ to $v_n$. 
\end{proposition}
{\bf Proof}
This follows from the fact that $\{\Phi^\pi: \pi\Vdash n\}$ is a basis of $\WSym_n$. Let $\varphi:\WSym_n\rightarrow V_n$ be the linear map such that $\varphi(\Phi^\pi)=V^\pi$; this is obviously a morphism of $\goth C$-module.  
$\Box$\\ \\
For convenience, we will write $\WSym[\varphi]=\varphi\WSym$.
\begin{example}\label{ExFree}
\rm
\begin{enumerate}
\item Let $A$ be an infinite alphabet. The restriction to $\WSym(A)$ of the morphism of algebra  $\varphi: \C\langle A\rangle\rightarrow \C[A]$ defined by $\varphi(a)=a$ is a morphism of $\goth C$-module sending $\WSym(A)$ to $Sym(A)$ (the algebra of symmetric functions on the alphabet $A$). 
This morphism can be defined without the help of alphabets, considering that $Sym$ is generated by the power sums $p_1,\dots, p_n,\dots$ with the action  $[\pi_1,\dots,\pi_m](p_{\#\pi_1},\dots,p_{\#\pi_m})=p_{\#\pi_1}\dots p_{\#\pi_m}$. 
\item Let $A$ be any alphabet (finite or not). If $V$ is a sub-$\goth C$-module of $\C\langle A\rangle$ generated by the homogeneous polynomials $P_n\in \C[A^n]$ as a $\goth C$-module, the linear map sending, for each $\pi=\{\pi_1,\dots,\pi_m\}$, $\Phi^\pi$ to $\shuffle_\Pi[P_{\#\pi_1},\dots,P_{\#\pi_m}]$, where $\Pi=[\pi_1,\dots,\pi_m]$, is a morphism of $\goth C$-module.
\end{enumerate}
\end{example}

\section{The Hopf algebra of set partitions into lists\label{BSym}}

\subsection{Set partitions into lists\label{SPL}}

A \emph{set partition into lists} is an object which can be constructed from a set partition by ordering each block. For example, $\{[1,2,3],[4,5]\}$ and $\{[3,1,2],[5,4]\}$ are two distinct set partitions into lists of the set $\{1,2,3,4,5\}$. The number of set partitions into lists of an $n$-element set (or set partitions into lists of size $n$) is given by Sloane's sequence A000262 \cite{Sloane}. If $\Pi$ is a set partition into lists of $\{1,\dots,n\}$, we will write $\Pi\Vvdash n$.
We will denote by $\mathrm{cycle\_support}(\sigma)$ the \emph{cycle support} of a permutation $\sigma$, \ie the set partition associated to its cycle decomposition. For instance, $\mathrm{cycle\_support}(325614)=\{\{135\},\{2\},\{4,6\}\}$. A set partition into lists can be encoded by a set partition and a permutation in view of the following easy result:
\begin{proposition}\label{prop:bw}
For all $n$, the set partitions into lists of size $n$ are in bijection with the pairs $(\sigma,\pi)$ where $\sigma$ is a permutation of size $n$ and $\pi$ is a set partition which is less fine than or equal to the cycle support of $\sigma$.
\end{proposition}
Indeed, from a set partition $\pi$ and a permutation $\sigma$, we obtain a set partition into lists $\Pi$ by ordering
the elements of each block of $\pi$ so that they appear in the same order as in $\sigma$.

\begin{example}\rm
Starting  the set partition $\pi=\{\{1, 4, 5\}, \{6\}, \{3, 7\}, \{2\}\}$ and the permutation $\sigma=4271563$, we obtain the set partition into lists $\Pi=\{[4, 1, 5], [7, 3], [6], [2]\}$.  
\end{example}

\subsection{Construction}

Let $\Pi\Vvdash n$ and $\Pi'\Vvdash n'$ be two set partitions into lists.
 Then, we set $\Pi \uplus \Pi' = \Pi \cup \{[l_1+n,\dots,l_k+n]:[l_1,\dots,l_k]\in\Pi'\}\Vvdash n+n'$.
   Let $\Pi'\subset\Pi\Vvdash n$, since the integers appearing in $\Pi'$ are all distinct, 
   the standardized $\std(\Pi')$ of $\Pi'$ is well defined as
    the unique set partition into lists obtained 
    by replacing the $i$th smallest integer in $\Pi$ by $i$. For example,
$\std(\{[5,2],[3,10],[6,8]\}) = \{[3,1],[2,6],[4,5]\}.$ 

\begin{definition}
The Hopf algebra $\BWSym$ is formally defined by its basis $(\Phi^\Pi)$ where the $\Pi$ are set partitions into lists, its product
$
\Phi^\Pi \Phi^{\Pi'} = \Phi^{\Pi\uplus \Pi'}
$
and its coproduct
$
\Delta(\Phi^\Pi) = \sum \Phi^{\std(\Pi')}\otimes\Phi^{\std(\Pi'')},
$
where the sum is over the $(\Pi',\Pi'')$ such that $\Pi'\cup \Pi''=\Pi$ and $\Pi'\cap\Pi'' = \emptyset$. 
\end{definition}
Following Section \ref{SPL} and for convenience, we will use alternatively $\Phi^\Pi$
 and $\Phi^{\binom\sigma \pi}$ to denote the same object.\\
We define $M_\Pi=M_{\binom\sigma \pi}$ by setting
$
\displaystyle \Phi^{\binom\sigma\pi}=\sum_{\pi\leq\pi'} M_{\binom\sigma\pi'}.
$
The formula being diagonal, it defines $M_\Pi$  for any $\Pi$ and proves that the family $(M_\Pi)_\Pi$ is a basis of $\BWSym$.
Consider for instance $ \{ [3,1],[2] \} \sim  \begin{pmatrix} 321 \\ \{\{1,3\},\{2\}\} \end{pmatrix} $,
we have 
$\Phi^{\{[3,1],[2]\}}= M_{\{[3,1],[2]\}} + M_{\{[3,2,1]\}}$.

For any set partition into lists $\Pi$, let $s(\Pi)$ be the corresponding classical set partition. Then, the linear application $\phi$ defined by 
$
\phi(\Phi^\Pi)=\Phi^{s(\Pi)}$
is obviously a morphism of Hopf algebras. As an associative algebra, $\BWSym$ has also algebraical links with the algebra $\FQSym$ \cite{FQSym}. Recall that this algebra is defined by its basis 
 $(E^\sigma)$ whose product is
$
E^\sigma E^\tau = E^{\sigma/\tau},
$
where $\sigma/\tau$ is the word obtained by concatening $\sigma$ and the word obtained from $\tau$ by adding the size of $\sigma$ to all the letters (for example $321/132=321564$).\\
The subspace $V$ of $\WSym \otimes \FQSym$ linearly spanned by the $\Phi^\pi \otimes E^\sigma$ such that the cycle supports of $\sigma$ is finer than $\pi$ is a subalgebra, and the linear application
sending $\Phi^{\pi,\sigma}$ to $\Phi_\pi \otimes E^\sigma$
is an isomorphism of algebras. 
Moreover, when the set of cycle supports of $\sigma$ is finer 
than $\pi$, $M_\pi \otimes E^\sigma$ also belongs to $V$ and is 
the image of $M_{\binom\sigma\pi}$.\\
The linear application from $\BWSym$ to $\FQSym$ which sends $M_{\{[\sigma]\}}$ to $E^\sigma$ and $M_\Pi$ to $0$ if ${\rm card}(\Pi)>1$, is also a morphism of algebras.

\subsection{Realization}
Let ${A}^{(j)}=\{a_i^{(j)} \mid i>0\}$ be an infinite set of bi-indexed noncommutative variables, with  
the order relation defined  by $a_i^{(j)} < a_{i'}^{(j)}$ if $i< i'$. Let $ A =\bigcup_j A ^{(j)}$. Consider the set partition into lists
$
\Pi = \{L^1,L^2,\ldots\} = \{[l^1_1, l^1_2,\ldots,l^1_{n_1}],[l^2_1,l^2_2,\ldots,l^2_{n_2}],\ldots\}\Vvdash n$. 
Then, one obtains a \emph{polynomial realization}  $\BWSym( A )$ by identifying $\Phi^\Pi$ with $\Phi^\Pi(A)$, 
the sum of all the monomials $a_1\ldots a_n$ (where the $a_i$ are in $ A $)
such that $k=k'$ implies  $a_{l^{k}_t},a_{l^{k'}_s}\in  A ^{(j)}$  for some $j$,  and for each $k$, 
$\std(a_{i_1}\dots a_{i_{n_k}})=\std(l^k_1\dots l^k_{n_k})$ with $\{l^k_1,\dots,l^k_{n_k}\}=\{i_1<\dots<i_{n_k}\} $
(The ``B'' of $\BWSym$ is for ``bi-indexed letters''). The coproduct $\Delta$ can be interpreted by identifying $\Delta(\Phi^\Pi)$ with $\Phi^\Pi ( A  +  B )$ as in the case of $\WSym$. Here, if $a\in A $ and $b\in  B $ then $a$ and $b$ are not comparable.
\begin{proposition}
The Hopf algebras $\BWSym$ and $\BWSym( A )$ are isomorphic.
\end{proposition}
Now, let $M'_\Pi(A)$
 be the sum of all the monomials $a_1\ldots a_n$, $a_i\in  A $, 
 such that   $a_{l^{k}_t}$ and $a_{l^{k'}_s}$ belong in the same $ A ^{(j)}$  if and only if $k=k'$,  and for each $k$, $\std(a_{i_1}\dots a_{i_{n_k}})=\std(l^k_1\dots l^k_{n_k})$ with $\{l^k_1,\dots,l^k_{n_k}\}=\{i_1<\dots<i_{n_k}\} $. For example, the monomial $a^{(1)}_1 a^{(1)}_1 a^{(1)}_2$ appears in the expansion of $\Phi_{\{[1,3],[2]\}}$, but not in the one of $M'_{\{[1,3],[2]\}}$. The $M'_\Pi$ form a new basis $(M'_\Pi)$ of $\BWSym$. Note that this basis is not the same as $(M_\Pi)$. For example, one has
\[
\Phi^{\{[1],[2]\}} = M_{\{[1],[2]\}}+M_{\{[1,2]\}} = M'_{\{[1],[2]\}}+M'_{\{[1,2]\}}+M'_{\{[2,1]\}}.
\]
Consider the basis $F_\sigma$ of $\FQSym$ defined in \cite{FQSym}. The linear application, from $\BWSym$ to $\FQSym$, which sends $M'_{\{[\sigma]\}}$ to $F_\sigma$, and $M'_\Pi$ to $0$ if ${\rm card} (\Pi)>1$, is a morphism of algebras.

\subsection{Related Hopf algebras}
By analogy with the construction of $\BWSym$, we define a ``B'' version for each of the algebras $\Sym$, $\QSym$ and $\WQSym$. In this section, we sketch briefly how to construct them; the complete study is deferred to a forthcoming paper.\\
As usual when $L=[\ell_1,\dots,\ell_k]$ and $M=[m_1,\dots,m_2]$ are
 two lists, the shuffle product is defined recursively by $[\,]\shuffle L=L\shuffle[\,]=\{L\}$ and $L\shuffle M=[\ell_1].([\ell_2,\dots,\ell_k]\shuffle M)\cup [m_1].(L\shuffle [m_2,\dots,m_k])$. 
The algebra of biword quasi-symmetric functions ($\BWQSym$) has its bases indexed by set compositions  into lists. 
The algebra is defined as the vector space spanned by the formal symbols 
$\Phi_\Pi$, where $\Pi$ is a composition into lists of the set $\{1,\dots,n\}$ for a given $n$, 
together with the product $\Phi_\Pi \Phi_{\Pi'}=\sum_{\Pi''\in\Pi\shuffle \Pi'[n] }\Phi_{\Pi''}$, 
where $\Pi'[n]$ means that we add $n$ to each of the integers in the lists
 of $\Pi'$ and $\Pi$ is a composition into lists of $\{1,\dots,n\}$. 
 Endowed with the coproduct defined by$\Delta(\Phi_{\Pi})=\sum_{\Pi'.\Pi''=\Pi}\Phi_{\std(\Pi')}\otimes \Phi_{\std(\Pi'')}$, $\BWQSym$ has a structure of Hopf algebra. 
 Note that $\BWQSym=\bigoplus_n\BWQSym_n$ is naturally graded; the dimension of the graded 
 component $\BWQSym_n$ is $2^{n-1}n!$ (see sequence A002866 in \cite{Sloane}).\\
The algebra $\BSym=\bigoplus_n\BSym_n$ is a graded cocommutative Hopf algebra whose bases are 
indexed by sets of permutations. 
Formally, we set $\BSym_n=\mathrm{span}\{\phi^{\{\sigma_1,\dots,\sigma_k\}}:\sigma_i\in\S_{n_i}, n_1+\dots+n_k=n\}$, 
$\phi^{S_1}.\phi^{S_2}=\phi^{S_1\cup S_2}$ and for any permutation $\sigma$, $\phi^{\{\sigma\}}$ is primitive. The dimensions of the graded components are given by the sequence A077365 of \cite{Sloane}.\\
Finally, $\BQSym=\bigoplus_n\BQSym_n$ is generated by 
$\phi_{[\sigma_1,\dots,\sigma_k]}$, its product is $\phi_L\phi_{L'}=\sum_{L''\in L\shuffle L'}\phi_{L''}$ 
and its coproduct is $\Delta(\phi_L)=\sum_{L=L'.L''}\phi_{L'}\otimes \phi_{L''}$. 
The dimension of the graded component $\BQSym_n$ is given by  Sloane's sequence A051296 \cite{Sloane}. 
\section{On the R-P theorem\label{Polyaetc}}

\subsection{R-P theorem and symmetric functions\label{PolyaSym}}

Consider two sets $X$ and $Y$ such that $X$ is finite ($\#X=n$), together with a group $G\subset \S_n$ acting on $X$ by permuting its elements.
We consider the set $Y^X$ of the maps $X\rightarrow Y$. The type of a map $f$ is the vector of the multiplicities of its images; more precisely, $\type(f)\in\N^Y$ with $\type(f)_y=\#\{x\in X:f(x)=y\}$. For instance, consider $X=\{a,b,c,d,e\}$, $Y=\{0,1,2\}$, $f(a)=f(c)=f(d)=1$, $f(b)=2$, $f(c)=0$: we have $\type(f)=[\mathop{1}^0,\mathop{3}^1,\mathop{1}^2]$.
The action of $G$ on $X$ induces an action of $G$ on $Y^X$. Obviously, the type of a function is invariant for the action of $G$. 
Then all the elements of an orbit of $G$ in $Y^X$ have the same type, so that the type of an orbit will be the type of its elements.
The question is: how to count the number ${\goth n}_I$ of orbits for the given type $I$?
 Note that, if $\lambda_I$ denotes the (integer) partition obtained by removing all the zeros in $I$ and reordering its elements in the decreasing order, $\lambda_I=\lambda_{I'}$ implies  ${\goth n}_I={\goth n}_{I'}$; it suffices to understand how to compute ${\goth n}_\lambda$ when $\lambda$ is a partition. The Redfield-P\'olya theorem deals with this problem and its main tool is the cycle index:
\[
Z_G:=\frac1{\#G}\sum_{\sigma\in G}p^{\mathrm{cycle\_type}(\sigma)},
\]
where $\mathrm{cycle\_type}(\sigma)$ is the (integer) partition associated to the cycle 
of $\sigma$ (for instance $\sigma=325614=(135)(46)$, $\mathrm{cycle\_type}(\sigma)=[3,2,1]$). When $\lambda=[\lambda_1,\dots,\lambda_k]$
 is a partition, $p^\lambda$ denotes the (commutative) symmetric function $p^\lambda=p_{\lambda_1}\dots p_{\lambda_k}$ and $p_n$ is the classical power sum symmetric function.\\ \\ \\
The Redfield-P\'olya theorem states:
\begin{theorem}\label{Polya}
The expansion of $Z_G$ on the basis  $(m_\lambda)$ of monomial symmetric functions is given by
\[
Z_G=\sum_\lambda {\goth n}_\lambda m_\lambda.
\]
\end{theorem} 
\begin{example}\label{polyav1}\rm
Suppose that we want to enumerate the non-isomorphic non-oriented graphs on three vertices.
 The symmetric group $\S_3$ acting on the vertices induces an action of the group 
 $$G:=\{123456,165432,345612,321654,561234,543216\}\subset \S_6$$ on the edges.
The construction is not unique. 
We obtain the group $G$ by labelling the $6$ edges from $1$ to $6$. 
Hence, to each permutation of the vertices corresponds a permutation of the edges. 
Here, the $1$ labels the loop from the vertex $1$ to itself, $2$ labels the edge which links the vertices $1$ and $2$, $3$ is the loop from the vertex $2$ to itself, $4$ labels the edge from the vertex $2$ to the vertex $3$, $5$ is the loop from the vertex $3$ to itself, finally, $6$ links the vertices $1$ and $3$.
The cycle index of $G$ is
\[\begin{array}{rcl}
Z_G&=&\frac16(p_1^6+3p_2^2p_1^2+2p_3^2)=m_6 + 2 m_{51} + 4 m_{42} + 6 m_{411} + 6 m_{33} + \dots
\end{array}
\]
The coefficient $4$ of $m_{42}$ means that there exists $4$ non-isomorphic graphs with $4$ edges coloured in blue and $2$ edges coloured in red.
\end{example}

\subsection{Word R-P theorem\label{PolyaWSym}}

If $\sigma$ is a permutation, we define
$
\Phi^\sigma:=\Phi^{\mathrm{cycle\_support}(\sigma)}.
$
 Now for our purpose, a map $f\in Y^X$ will be encoded by a word $w$ : 
 we consider  an alphabet $A=\{a_y:y\in Y\}$, the elements of $X$ are relabelled by $1,2,\dots,\# X=n$
  and $w$ is defined as the word $b_1\dots b_n\in A^n$ such that $b_i=a_{f(i)}$. With these notations, the action of $G$ on $Y^X$ is encoded by the action of the permutations of $G$ on the positions of the letters in the words of $A^n$.\\
It follows that for any permutation $\sigma\in G$, one has
\begin{equation}\label{Invariance1}
\Phi^\sigma=\sum_{w\sigma=w}w.
\end{equation}
The cycle support polynomial is defined by
$
{\mathrm Z}_G:=\sum_{\sigma\in G} \Phi^\sigma. 
$
From (\ref{Invariance1}) we deduce 
\[
{\mathrm Z}_G=\sum_{w}\#\Stab_G(w)w
\]
where $\Stab_G(w)=\{\sigma\in G:w\sigma=w\}$ is the subgroup of $G$ which stabilizes $w$.
In terms of monomial functions, this yields :
\begin{theorem}\label{WPolya}
\[{\mathrm Z}_G=\sum_{\pi}\#\mathrm{Stab}_G(w_\pi)M_\pi \]
where $w_\pi$ is any word $a_1\dots a_n$ such that $a_i=a_j$ if and only if $i,j\in\pi_k$ for some $1\leq k\leq n$.
\end{theorem}

\begin{example}\label{polyav2}\rm
Consider the same example as in Example \ref{polyav1}. Each graph is now encoded by a word $a_1a_2a_3a_4a_5a_6$: the letter $a_1$ corresponds to the colour of the vertex $1$, the letter $2$ to the colour of the vertex $2$ and so on.\\
The cycle support polynomial is
{\footnotesize\[\begin{array}{rcl}
{\mathrm Z}_G&:=&\Phi^{\{\{1\},\{2\},\{3\},\{4\},\{5\},\{6\}\}}+\Phi^{\{\{2,6\},\{3,5\},\{1\},\{4\}\}}
\\&&+\Phi^{\{\{1,3\},\{4,6\},\{2\},\{5\}\}}+\Phi^{\{\{1,5\},\{2,4\},\{3\},\{6\}\}}
+2\Phi^{\{\{1,3,5\},\{2,4,6\}\}}.
\end{array}\]}
The coefficient of $M_{\{\{2,6\},\{3,5\},\{1\},\{4\}\}}$ in ${\mathrm Z}_G$ is $2$ because it
 appears only in $\Phi^{\{\{1\},\{2\},\{3\},\{4\},\{5\},\{6\}\}}$ and $\Phi^{\{\{2,6\},\{3,5\},\{1\},\{4\}\}}$.
  The monomials of $M_{\{\{2,6\},\{3,5\},\{1\},\{4\}\}}$ are of the form $abcdba$, where $a, b, c$ and $d$ are four distinct letters. The stabilizer of $abcdba$ in $G$ is the two-element subgroup $\{123456,165432\}$.
\end{example}

\subsection{From word R-P theorem to R-P theorem}

The aim of this section is to link the numbers ${\goth n}_I$ of  Section \ref{PolyaSym} and the numbers $\#\Stab_G(w)$ appearing in Section \ref{PolyaWSym}.\\
If $w$ is a word we will denote by $\orb_G(w)$ its orbit under the action of $G$. The Orbit-stabilizer theorem (see \emph{e.g.}\cite{Burnside}) together with  Lagrange's theorem gives:
\begin{equation}\label{Burnside}
 \#G=\#\orb_G(w)\#\Stab_G(w)
\end{equation}
Denote by $\Lambda(\pi)$ the unique integer partition defined by $(\#\pi_1,\dots,\#\pi_k)$ if $\pi=\{\#\pi_1,\dots,\pi_k\}$ with $\#\pi_1\geq\#\pi_2\geq\dots\geq\#\pi_k$. If $\lambda=(m^{k_m},\dots,2^{k_2},1^{k_1})$ we set $\lambda^!=k_m!\dots k_2!k_1!$.
The shape of a word $w$ is the unique set partition $\pi(w)$ such that $w$ is a monomial of $M_{\pi(w)}$.
 Note that all the orbits of words with a fixed shape $\pi$ have the same cardinality.
  Furthermore, let $\pi=\{\pi_1,\dots,\pi_k\}$ and $A_k$ be an alphabet of size $k$, 
  the number of words of shape $\pi$ on $A_k$ equals $\Lambda(\pi)^!$. Hence the set of all the words of shape $\pi$ on the alphabet $A_k$ is partitioned into $\frac{\Lambda(\pi)^!}{\#\orb(w_\pi)}$ orbits of size $\#\orb(w_\pi)$. We deduce that:
\begin{equation}\label{ntoStab}
{\goth n}_\lambda=\sum_{\Lambda(\pi)=\lambda}\frac{\lambda^!}{\#\orb_G(w_\pi)}=\sum_{\Lambda(\pi)=\lambda}\frac{\lambda^!\#\Stab_G(w_\pi)}{\#G}.
\end{equation}
If we consider the morphism of algebra $\theta:\WSym\rightarrow Sym$ which sends $\Phi_n$ to $p_n$,
 we have $\theta(M_\pi)=\Lambda(\pi)^{!}m_\lambda$. Hence, we have
\[
 \frac1{\#G}\theta({\mathrm Z}_G)=\sum_\lambda \left(\sum_{\Lambda(\pi)=\lambda}\frac{\lambda^!}{ \#\orb_G(w_\pi)}\right)m_\lambda=\sum_\lambda {\goth n}_\lambda m_\lambda
\]
as expected by the Redfield-P\'olya theorem (Theorem \ref{Polya}).

\subsection{R-P theorem without multiplicities}
In this section we give a $\BWSym$ version of the  R-P theorem whose main property is to have no multiplicities.
Examining with more details Example \ref{polyav2},  
the coefficient $2$ of $M_{\{2,6\},\{3,5\},\{1\},\{4\}}$ in $\mathrm Z_G$ follows
 from the  group $\{123456,165432\}$ of order two which stabilizes $abcdcb$.
In terms of set partitions into lists, this can be interpreted by 
$
\ M_{\{[2,6],[3,5],[1],[4]\}} + M_{\{[6,2],[5,3],[1],[4]\}}\rightarrow 2M_{\{\{2,6\},\{3,5\},\{1\},\{4\}\}}.
$
 We deduce the following version (without multiplicities) of Theorem \ref{WPolya} in $\BWSym$.
\begin{theorem}\label{Polyawithoutm}
Let $G$ be  a permutation group. We have $$\mathbb Z_G:=\sum_{\sigma \in G}  \Phi^{ \binom{\sigma}{ 
\mathrm{cycle}\_\mathrm{support}(\sigma)}}=\sum_{\pi} \sum_{\sigma \in Stab_{\pi}(G)}  M_{ \binom{\sigma}\pi}.$$
\end{theorem} 
Consider again Example \ref{polyav2}. 
{
\[\begin{array}{rl}
\mathbb Z_{G}= & \Phi^{  \binom{123456}{\{ \{1\}\{2\}\{3\}\{4\}\{5\}\{6\} \} } } +
 \Phi^{  \binom{165432}{\{ \{1\}\{26\}\{35\}\{4\} \}} } + \Phi^{\binom{345612} {\{ \{135\}\{246\} \}} }
\\+ & \Phi^{  \binom{321654}{ \{ \{13\}\{2\}\{46\}\{5\} \}} } + 
\Phi^{ \binom{ 561234}{\{ \{135\}\{246\} \}} } + 
\Phi^{  \binom{543216}{\{ \{15\}\{234\}\{6\} \}}}.\end{array}\] }
When expanded in the monomial $M$ basis, there are exactly $2$ terms of the form $M_{\binom\sigma{\{\{2,6\},\{3,5\},\{1\},\{4\}\}}}$ (for $\sigma=123456$ and $\sigma=165432$). 
Note that we can use another realization which is compatible with the space but not with the Hopf algebra structure.
It consists to set $\widetilde\Phi^{\binom\sigma\pi}:=\sum_w\binom\sigma w$, 
where the sum is over the words $w = w_1\dots w_n$ ($w_i \in  A $) such that if $i$ and $j$
are in the same block of $\pi$ then $w_i = w_j$.
If we consider the linear application $\widetilde\psi$ sending $\Phi^{\binom\sigma\pi}$ 
to $\widetilde\Phi^{\binom\sigma\pi}$, $\widetilde\psi$ sends $M_{\binom\sigma\pi}$ to
$\widetilde M_{\binom\sigma\pi}:=\sum_w\binom\sigma w$, 
where the sum is over the words $w = w_1\dots w_n$ ($w_i \in  A $) such that  $i$ and $j$
are in the same block of $\pi$ if and only if $w_i = w_j$. 
Let $w$ be a word, the set of permutations $\sigma$ such that $\binom\sigma w$ appears in the expansion of 
$\widetilde\psi(\mathbb Z_G)$ is the stabilizer of $w$ in $G$.
The linear application sending each
 biword $\binom\sigma w$ to $w$ sends $\widetilde \Phi^{\binom\sigma \pi}$ to $\Phi^\pi$ and 
 $\sum_{\sigma\in\Stab_G(w)}\binom\sigma w$ to $\# \Stab_G(w)w$. Note that $\# \Stab_G(w)$ is also the coefficient
 of the corresponding monomial $M_{\pi(w)}$ in the cycle support polynomial $\mathrm Z_G$. 
 For instance, we recover the coefficient $2$ in 
 Example \ref{polyav2} from the biwords $\binom{123456}{abcdcb}$
 and $\binom{165432}{abcdcb}$ in $\widetilde\psi(\mathbb Z_G)$.

\subsection{$\WSym$ and Harary-Palmer type enumerations\label{Harary}}

Let $ A :=\{a_1,\dots,a_m\}$ be a set of formal letters and $I=[i_1,\dots,i_k]$ a sequence of elements of $\{1,\dots,m\}$. We define the virtual alphabet $ A _{I}$ by
\[
\Phi^{n}( A _{I}):=(a_{i_1}\dots a_{i_k})^{\frac n k}+(a_{i_2}\dots a_{i_k}a_{i_1})^{\frac n k}+\dots
+(a_{i_{k}}a_1\dots a_{i_{k-1}})^{\frac n k},
\]
if $k$ divides $n$ and $0$ otherwise. If $\sigma\in\S_m$ we define the alphabet $ A _\sigma$ as the formal sum of the alphabets $ A _c$ associated to its cycles:
\[
\Phi^{\{1\dots n\}}[ A _\sigma]:=\sum_{c\mbox{ cycle in }\sigma}\Phi^{n}[ A _c].
\]
From Example \ref{ExFree}.2, the set $\{\Phi^{\{1,\dots,n\}}[ A _\sigma]:n\in\N\}$ generates the sub-$\goth C$-module $\WSym[ A _\sigma]$ of $\C\langle  A \rangle$ (the composition $\Pi$ acting by $\shuffle_\Pi$).\\
Let $H\subset \S_m$ and $G\subset \S_n$ be two permutation groups. We define $
{\mathrm Z}(H;G):=\sum_{\tau\in H}\Phi^G[ A _\tau].
$

\begin{proposition}

We have:
$${\mathrm Z}(H;G)=\sum_{w\in A ^n}\#{\mathrm Stab}_{H,G}(w)w$$
where ${\mathrm Stab}_{H,G}(w)$ denotes the stabilizer of $w$ under the action of $H\times G$ ($H$ acting on the left on the names of the variables $a_i$ and $G$ acting on the right on the positions of the letters in the word); equivalently, 
${\mathrm Stab}_{H,G}(a_{i_1}\dots a_{i_k})=\{(\tau,\sigma)\in H\times G:a_{\tau(i_{\sigma(j)})}\mbox{ for each }1\leq j\leq n\}$. 
\end{proposition}
Hence, from Burnside's classes formula, sending each variable to $1$ in ${\mathrm Z}(H;G)$, we obtain the number of orbits of $H\times G$.
\begin{example}\rm
Consider the set of the non-oriented graphs without loop whose edges are labelled by three colours. Suppose that we consider the action of the group $H=\{123,231,312\}\subset\S_n$ on the colours. We want to count the number of graphs up to permutation  of the vertices ($G=\S_3$) and the action of $H$ on the edges. There are three edges, and each graph will be encoded by a word $a_{i_1}a_{i_2}a_{i_3}$ where ${i_j}$  denotes the colour of the edge $j$. We first compute the specialization $\Phi^{\{1\dots n\}}[ A _\sigma]$  for $1\leq n\leq 3$ and $\sigma\in H$. We find $\Phi^{\{1\}}[ A _{123}]=a_1+a_2+a_3$,  $\Phi^{\{1\}}[ A _{231}]=\Phi^{\{1\}}[ A _{312}]=0$, $\Phi^{\{1,2\}}[ A _{123}]=a_1^2+a_2^2+a_3^2$,  $\Phi^{\{1,2\}}[ A _{231}]=\Phi^{\{1,2\}}[ A _{312}]=0$,  $\Phi^{\{1,2,3\}}[ A _{123}]=a_1^3+a_2^3+a_3^3$, $\Phi^{\{1,2,3\}}[ A _{213}]=a_2a_3a_1+a_3a_1a_2+a_1a_2a_3$, and
$\Phi^{\{1,2,3\}}[ A _{312}]=a_1a_3a_2+a_3a_2a_1+a_2a_1a_3$. Now, we deduce the values of the other $\Phi^{\pi}[ A _{\sigma}]$ with $\pi\Vdash 3$ and $\sigma\in H$ by the action of $\shuffle_\Pi$. For instance:
$$
\begin{array}{rcl}
\Phi^{\{\{1,2\},\{3\}\}}[ A _{123}]&=&\shuffle_{[\{1,2\},\{3\}]}\left(\Phi^{\{1,2\}}[ A _{123}],\Phi^{\{1\}}[ A _{1}]\right)\\&=&a_1^3+a_1^2a_2+ a^2_1a_3 + a^2_2 a_1 + a^3_2 + a^2_2a_3 + a^2_3a_1 + a^2_3 a_2 + a^3_3.
\end{array}$$
We find also
\begin{align*}
\Phi^{\{\{1,3\},\{2\}\}}[ A _{123}] &= a^3_1 +a_1a_2a_1+a_1a_3a_1+a_2a_1a_2+a^3_2 +a_2a_3a_2+a_3a_1a_3+a_3a_2a_3+a^3_3,\\
\Phi^{\{\{1\},\{2,3\}\}}[ A _{123}] &= a^3_1 + a_2a^2_1 + a_3a^2_1 + a_1a^2_2 + a^3_2 + a_3a^2_2 + a_1a^2_3 + a_2a^2_3 + a^3
_3,\\
\Phi^{\{\{1\},\{2\},\{3\}\}}[ A _{123}] &= (a_1 + a_2 + a_3)^3.
\end{align*}
The other $\Phi^\pi[ A _\sigma]$ are zero. Hence,
\[\begin{array}{rcl}
\mathrm Z[H;\S_3]&=&
\Phi^{123}[ A _{123}] + \Phi^{132}[ A _{123}] +
 \Phi^{213}[ A _{123}] + \Phi^{321}[ A _{123}]+\\&&
\Phi^{231}[ A _{231}] + \Phi^{231}[ A _{312}] + \Phi^{312}[ A _{231}] 
+ \Phi^{312}[ A _{312}]\\
&=&6(a^3_1 + a^3_2 + a^3_3) + 2\sum_{
i\neq j} a^2_ia_j + 2
\sum_{i\neq j} a_ja_ia_j + 
2\sum_{i\neq j} a_ja^2_i\\&&
+3(a_1a_2a_3 + a_2a_3a_1 + a_3a_1a_2) + 3(a_1a_3a_2 + a_3a_2a_1 + a_2a_1a_3)
.
\end{array}
\]
The coefficient $3$ of $a_1a_2a_3$ means that the word is invariant under the action of three pairs of permutations (here $(123,123),\ (231,312),\ (312,231)$).
 Setting $a_1=a_2=a_3=1$, we obtains $\mathrm Z[H,\S_3]=18\times 4$: $18$ is the order of the group $H\times \S_3$ and $4$ is the number of orbits:\\
$\{a^3_1, a^3_2, a^3_3\}$,
$\{a^2_1a_2, a_1a_2a_1, a_2a^2_1, a^2_2a_3, a_2a_3a_2, a_3a^2_2, a^2_3a_1, a_3a_1a_3, a_1a^2
_3\}$,
$\{a^2_2a_1, a_2a_1a_2, a_1a^2_2, a^2_1a_3$, $ a_1a_3a_1,a_3a^2_1, a^2_3a_2, a_3a_2a_3, a_1a^2_3\}$, and 
$\{a_1a_2a_3, a_1a_3a_2, a_2a_1a_3, a_2a_3a_1, a_3a_1a_2, a_3a_2a_1\}$.
 \end{example}

\end{document}